\documentclass[10pt,final]{amsart}
\usepackage{amsmath,amssymb,amsfonts,amsthm,slashed,braket}
\usepackage[all]{xy}
\usepackage{mathrsfs}
\usepackage{showkeys}
\usepackage{hyperref}

\title{Further Counterexamples to the Integral Hodge Conjecture}
\author{Arnav Tripathy}

\newcommand\nc{\newcommand}
\nc\linesep{\bigskip}
\nc\newprob[1]{\marginnote{#1}[\parskip]}
\nc\bA{\mathbb A}
\nc\bC{\mathbb C}
\nc\bD{\mathbb D}
\nc\bR{\mathbb R}
\nc\bZ{\mathbb Z}
\nc\bQ{\mathbb Q}
\nc\bP{\mathbb P}
\nc\bV{\mathbb V}
\nc\bW{\mathbb W}
\nc\bG{\mathbb G}
\nc\mf\mathfrak
\nc\mc\mathcal
\nc\mb\mathbb
\nc\brac[1]{\langle#1\rangle}
\nc\abs[1]{\lvert#1\rvert}
\nc\norm[1]{\lVert#1\rVert}
\nc\onto{\twoheadrightarrow}
\nc\into{\hookrightarrow}
\nc\lto{\longrightarrow}
\nc\action{\curvearrowright}
\nc\on\operatorname
\nc\wbar\overline
\nc\what\widehat
\nc\wtilde\widetilde
\nc\nop\DeclareMathOperator
\nc\eps{\varepsilon}
\nc\tsym{\widetilde{\text{Sym}}}
\nc\oarrow[1]{\overset{#1}\to}
\nop\Hom{Hom}
\nop\End{End}
\nop\Aut{Aut}
\nop\im{Im}
\nop\id{id}
\nop\tr{Tr}
\nop\coker{coker}
\nop\Spec{Spec}
\nop\Jac{Jac}
\nop\Ext{Ext}
\nop\Tor{Tor}
\nc\op{\text{op}}
\nop\loc{Loc}
\nop\Frac{Frac}
\nc\ann{\text{ann}}
\nop\QCoh{QCoh}
\nop\Coh{Coh}
\nop\Sym{Sym}
\nop\gr{Gr}
\nop\Tot{Tot}
\nop\Fl{Fl}
\nop\tGamma{\widetilde\Gamma}
\nop\tloc{\widetilde{\text{Loc}}}
\nop\rep{Rep}
\nop\proj{Proj}
\nc\oo[1]{\overset\circ{#1}}
\nop\ospec{\oo{Spec}}
\nop\oTot{\oo{Tot}}
\nop\Bl{Bl}
\nop\Comp{Comp}
\nop\Ho{Ho}
\nop\cone{Cone}
\nop\LKE{LKE}
\nop\RKE{RKE}
\nop\pd{pd}
\nop\cd{cd}
\nop\depth{depth}
\nop\ass{Ass}
\nop\supp{supp}
\nop\codim{codim}
\nop\holim{\underset{\lto}{holim}}
\nop\dlim{\underset{\lto}{lim}}

\nop\uHom{\underline{\Hom}}
\nop\Pic{Pic}
\nop\Cl{Cl}
\nop\Div{Div}
\nop\rank{rank}
\nop\Der{Der}
\nop\dimrel{dim.rel}
\nc\sHom{\mathscr Hom}
\nc\sExt{\mathscr Ext}
\nc\dto{\dashrightarrow}
\nop\rspec{\bf Spec}
\nop\Gal{Gal}
\nop\Ind{Ind}
\nop\Frob{Frob}
\nop\Fib{Fib}
\nop\ratdim{rat\ dim}
\nop\Mod{Mod}
\nop\rat{rat}
\nop\val{val}
\nop\Rep{Rep}
\nop\colim{colim}

\theoremstyle{theorem}
\newtheorem{thm}{Theorem}

\theoremstyle{remark}

\begin{document} 
\maketitle
\tableofcontents

\begin{abstract}

We study the integral Hodge conjecture in complex codimension $2$ and $3$ for approximations to the classifying space of groups of type A. In degree two, we prove a conjecture of Ben Antieau, extending his two counterexamples to a general family of groups. In particular, when reduced to a finite field, these spaces extend Antieau's counterexamples to the integral $\ell$-adic Tate conjecture from the cases $\ell = 2, 3$ to all primes $\ell$.

\end{abstract}

\section{Introduction}
We continue the esteemed tradition of providing counterexamples to the integral Hodge conjecture that for $X$ a smooth, projective complex variety, the natural cycle class map $CH^k(X) \to H^{2k}(X; \mb{Z}) \cap H^{k,k}(X)$ is an isomorphism, where as usual by the target above, we mean more precisely the preimage of the Hodge summand $H^{k,k}(X)$ under the natural map $H^{2k}(X; \mb{Z}) \to H^{2k}(X; \mb{C})$. The first counterexamples were provided by Atiyah and Hirzebruch in~\cite{Atiyah-Hirzebruch} and were re-interpreted and extended in~\cite{TotaroCycleClass} by Totaro's refined cycle class map $CH^*(X) \to MU^*(X) \otimes_{MU_*} \mb{Z}$, where $MU$ is the extraordinary cohomology theory of complex cobordism; these counterexamples were torsion cycles on algebraic approximations to classifying spaces of algebraic groups. Kollar provided the first non-torsion counterexamples in~\cite{KollarHodge} using the quite different method of degenerations. In this paper, we hew very closely to the first line of thinking, especially as spurred on by the recent papers~\cite{Pirutka-Yagita},~\cite{Kameko}, and~\cite{Antieau} that extend further counterexamples to the integral Hodge conjecture to also be counterexamples to the integral Tate conjecture over finite fields. In particular, the current paper is almost entirely inspired by~\cite{Antieau}, in which Antieau provides new counterexamples to the integral Hodge and Tate conjectures that are strikingly different from the original Atiyah-Hirzebruch examples, and all successive examples built from approximations to classifying spaces. Indeed, in all previous counterexamples to the $\ell$-adic Tate conjecture, the Atiyah-Hirzebruch spectral sequence differential $d_{2\ell-1}$ does not vanish, essentially because these counterexamples were built in a very similar way from the original Atiyah-Hirzebruch counterexamples. In this paper, we prove Conjecture $1.3$ of~\cite{Antieau}, which we reproduce below:

\begin{thm} For $G = SL_{\ell^2} / \mu_{\ell}$, where $\ell$ is an odd prime, the image of $CH^2(BG) \to H^4(BG; \mb{Z}) \simeq \mb{Z} c_2$ is generated by $\ell c_2$. \end{thm}

Establishing this conjecture allows an extension of Antieau's counterexamples from the primes $\ell = 2, 3$, which he established by direct, exhaustive computation, to all primes $\ell$, using a more conceptual technique ideally well-suited for further investigation of counterexamples. In particular, we provide further counterexamples to the integral Hodge conjecture for $CH^3$ mapping to $H^6(-, \mb{Z})$ of (approximations to) classifying spaces and outline an approach for attempting to find counterexamples to Yagita's conjecture $12.2$ of~\cite{Yagita} that $CH^*(BG) \to MU^*(BG) \otimes_{MU_*} \mb{Z}$ is an isomorphism. We hope to return to this question and the interesting structure of the differentials in the Atiyah-Hirzebruch spectral sequence on these counterexamples using our techniques in future work.

The author learned of this work from conversations with Ben Antieau at the 2015 Utah Algebraic Geometry Summer Research Institute and would like to thank the organizers of the institute, as well as the preceding Bootcamp, for providing a wonderful and highly informative environment. Particular thanks are, of course, due to Ben Antieau for generously explaining his work on the problem and discussing several further stimulating ideas along similar lines. Antieau then further informed me later that Masaki Kameko also communicated a similar solution to Antieau's original conjecture, so do also see Kameko's paper for related ideas and approaches.

\section{Background and Notation}

We will be concerned in this paper with groups of type $A$, of the form $G = SL_{n} / \mu_{\ell}$, for $\ell$ a prime and $n = k \ell$ a multiple of $\ell$. Certainly, the techniques we use can be extended to more general cases, such as $\ell$ not a prime, for example, but we focus on this relatively restricted case. We will also consider the map $CH^*(BG) \to H^*(BG; \mb{Z})$, where the source of the map is Totaro's Chow groups of classifying spaces, as developed in~\cite{TotaroClassifying}. Moreover, we will typically be concerned only with the quotients of the above groups that survive under pullback to $BSL_n$; in other words, we are principally interested in the images of the top row in the bottom row under the following commutative square: \begin{equation*} \xymatrix{ CH^*(BG) \ar[d] \ar[r] & H^{2*}(BG; \mb{Z}) \ar[d] \\ CH^*(BSL_n) \ar[r]^{\sim} & H^{2*}(BSL_n; \mb{Z}) }. \end{equation*} In particular, whenever we consider the image of some group, we mean the image under the composition of the relevant cycle class map and pullback to the integral cohomology of $BSL_n$.

We first start by considering $* = 2$ above. The easier question to address is the image of $H^4(BG; \mb{Z})$ inside $H^4(BSL_n; \mb{Z}) \simeq \mb{Z} c_2$. The first case to understand is where $k = 1$, i.e. the case of $H^4(BPGL_{\ell}; \mb{Z}) \to H^4(BSL_{\ell}; \mb{Z})$. Following, for example, Lemma $3.2$ of~\cite{Antieau-Williams}, a straightforward Chern class computation yields that the image has index $\ell$ for $\ell$ odd and index $4$ for $\ell = 2$; for simplicity, we denote this number by $\ell'$ rather than some contrived expression like $(1 + \delta_{2, \ell})\ell$. Now, for general $k$, we consider the Serre spectral sequence for $BSL_n \to BG \to K(\mb{Z}/\ell, 2)$; $H^4(BG; \mb{Z})$ is identified with the kernel of the differential $d_4: H^4(BSL_n, \mb{Z}) \to H^5(K(\mb{Z}, \ell), \mb{Z})$. To understand this transgression, we consider the map $K(\mb{Z}/\ell, 1) \to BSL_n$. The diagonal embedding $\mu_{\ell} \to SL_n$ may be considered as a block-diagonal embedding of $\mu_{\ell}$ into $k$ successive blocks of $SL_{\ell}$. As each inclusion of these subgroups $SL_{\ell} \hookrightarrow SL_n$ are homotopic to each other, the total induced map $K(\mb{Z}/\ell, 1) \to BSL_n$ is equal in the homotopy category to $k$ times the composite $K(\mb{Z}/\ell, 1) \to BSL_{\ell} \to BSL_n$. In other words, $d_4$ for a general $k$ is $k$ times the composite $H^4(BSL_n, \mb{Z}) \stackrel{\sim}{\to} H^4(BSL_{\ell}, \mb{Z}) \stackrel{d_4}{\to} H^5(K(\mb{Z}/\ell, 2), \mb{Z})$, where the latter map is the differential $d_4$ for the case $k = 1$. We hence obtain that $\ker d_4 \subset H^4(BSL_n, \mb{Z})$ is generated by $\ell' / \gcd(k, \ell')$. Indeed, this computation motivated Antieau in his paper to look for examples with $k$ divisible by $\ell'$ (in fact, he took the first case of $k = \ell'$) as then the pullback $H^4(BG; \mb{Z}) \to H^4(BSL_n, \mb{Z})$ is an isomorphism. After all, if we wish to find counterexamples to the integral Hodge conjecture, our plan should be to attempt to maximize the index of $\im CH^2(BG)$ while minimizing the index of $\im H^4(BG; \mb{Z})$. 

We now turn to the question of determining the image of $CH^2(BG)$. In fact, the cycle class map is injective, by Lemma $15.1$ of~\cite{TotaroBook}, but we are more concerned with the fact that it is generated by Chern classes of representations, as explained in the same lemma. Indeed, in general, we will find it useful to define the subgroup $Ch^*(BG) \subset CH^*(BG)$ generated by Chern classes of complex representations: to be clear, a complex representation determines a complex vector bundle on (the approximations to) $BG$ and hence we may consider its Chern classes inside the Chow ring. In general, Riemann-Roch without denominators, as explained in Theorem $2.25$ of~\cite{TotaroBook}, yields that we have $(i-1)! CH^i(BG) \subset Ch^i(BG) \subset CH^i(BG)$, which we will use more seriously in the next section when we turn to $i = 3$. Here, we simply accept that $Ch^2(BG) \simeq CH^2(BG)$ and turn entirely to considerating the representation theory of $G$.

We attempt to establish some notational clarity in order to ideally reduce some awkward confluences at this point. Given a representation $\rho$ of $G$, we will be interested in its Chern classes, typically as pulled back to the integral cohomology of $BSL_n$. As such, these Chern classes will be expressed in terms of the cohomology of $BSL_n$, which is also written in terms of Chern classes. To attempt to minimize confusion, we will denote the latter Chern classes by the typical notation of lowercase $c_i$s, i.e. we write $H^*(BSL_n, \mb{Z}) = \mb{Z}[c_2, \cdots, c_n]$ as usual. We will denote Chern classes of representations by capital $C_i$s, and so we will typically be writing equations expressing $C_i(\rho)$ in terms of the $c_i$s. Finally, we reiterate the definition of $C_i(\rho)$: given a homomorphism $\rho: G \to GL_N$, we can deloop to obtain a morphism $BG \to BGL_N$; we then consider the composition $$C_i(\rho): BSL_n \to BG \to BGL_N \stackrel{c_i}{\to} K(\mb{Z}, 2i).$$ After carefully attempting to disambiguate our notation, we will immediately introduce an abuse thereof: as $H^4(BSL_n; \mb{Z}) \simeq \mb{Z} c_2$ and $H^6(BSL_n; \mb{Z}) \simeq \mb{Z} c_3$, we will equivocate between designating by $C_2(\rho)$ the actual fourth-degree cohomology class, as we defined it above, and instead simply the integer such that the cohomology class is that multiple of $c_2$. In later sections, we will do the same for $C_3$.

Recall that the representation ring of $SL_n$ is a truncated ring of symmetric functions, and we will frequently conflate the two, denoting representations by their associated symmetric polynomials. This ring has an integral basis indexed by Young tableaux with fewer than $n$ rows. The representation ring of $G$ will be the subring spanned by representations which are trivial on the diagonal $\mu_{\ell}$; in terms of the Young tableau, this condition translates to the total number of boxes being divisible by $\ell$, and so $\Rep(G)$ has an integral basis indexed by this subset of Young tableaux. There are a few alternative approaches to exactly how the Young tableaux parametrize representations. We find it convenient here to work with the elementary basis, as opposed to the somewhat more standard Schur basis. More explicitly, we describe $\Rep(SL_n) \simeq \mb{Z}[e_1, \cdots, e_{n-1}]$, the integral polynomial algebra on the antisymmetric polynomials of degree less than $n$. The subring $\Rep(G)$ now has an integral basis given by monomials $\prod_{\alpha \in S} e_{\alpha}$, where $S$ is a multiset of possibly repeated indices between $1$ and $n-1$, subject to the condition that $\sum_{\alpha \in S} \alpha$ is divisible by $\ell$.

The Chern character is our favored tool for computations of the $C_i$s. Notably, the Chern character, which we consider as the composition $ch: \Rep(G) \to K(BG) \to H^*(BG; \mb{Q}) \to H^*(BSL_n, \mb{Q})$, is a ring homomorphism, where the first morphism is the Atiyah-Segal completion morphism and the second is the usual Chern character. We now recall the expression of the Chern character in terms of Chern classes, which is given by writing the symmetric function $\sum \exp(t_a)$ in terms of the elementary totally antisymmetric polynomials of the $t_a$; for example, the first few terms are as follows: \begin{eqnarray*} \sum \exp(t_a) &=& \rank + \Big( \sum t_a \Big) + \Big( \frac{1}{2} ( \sum t_a )^2 - \sum_{a < b} t_a t_b \Big) + \\ && \Big( \frac{1}{6} (\sum t_a)^3 - \frac{1}{2} (\sum t_a) (\sum_{a < b} t_a t_b) + \frac{1}{2} \sum_{a < b < c} t_at_bt_c \Big) + O(t^4). \end{eqnarray*} Applying this to a representation $\rho$ as above, and noting that the first Chern class vanishes as we are working with the cohomology of $SL_n$ rather than $GL_n$, we have \begin{equation} \label{chern} ch(\rho) = \rank(\rho) - C_2(\rho) + \frac{1}{2} C_3(\rho) + \cdots. \end{equation}

\section{Second Chern Classes}

We turn in earnest to computing second Chern classes of representations, starting with the totally antisymmetric $e_i$s. We use the splitting principle, so let $X$ be a space mapping to $BSL_n$ such that the pullback on cohomology is injective and the standard representation splits into the sum of line bundles $\mc{L}_1 \oplus \cdots \oplus \mc{L}_n$. Then the composite $X \to BSL_n \stackrel{e_i}{\to} BSL_N$ is given by the sum of line bundles $$\bigoplus_{\lambda_1 < \cdots < \lambda_i} \mc{L}_{\lambda_1} \otimes \cdots \otimes \mc{L}_{\lambda_i}.$$ If we denote $c_1(\mc{L}_a) = t_a$, then we have \begin{eqnarray*} ch(e_i) &=& \sum_{\lambda_1 < \cdots < \lambda_i} ch(\mc{L}_{\lambda_1}) \cdots ch(\mc{L}_{\lambda_i}) \\ &=& \sum_{\lambda_1 < \cdots < \lambda_i} \exp(t_{\lambda_1}) \cdots \exp(t_{\lambda_i}) \\ &=& \sum_{\lambda_1 < \cdots < \lambda_i} \exp(t_{\lambda_1} + \cdots + t_{\lambda_i}) \\ &=& \rank(e_i) + \Big( \binom{n-1}{i-1} \sum t_a \Big) + \Big( \frac{1}{2} \binom{n-1}{i-1} \sum t_a^2 + \binom{n-2}{i-2} \sum_{a < b} t_a t_b \Big) + \\ && \Big( \frac{1}{6} \binom{n-1}{i-1} \sum t_a^3 + \frac{1}{2} \binom{n-2}{i-2} \sum_{a \ne b} t_a^2 t_b + \binom{n-3}{i-3} \sum_{a < b < c} t_a t_b t_c \Big) + O(t^4). \end{eqnarray*} We use~\ref{chern} to identify the second component: \begin{equation} \label{2} -C_2(e_i) = \frac{1}{2} \binom{n-1}{i-1} (\sum t_a)^2 + \Big( \binom{n-2}{i-2}  - \binom{n-1}{i-1} \Big) \sum_{a < b} t_at_b, \end{equation} whereupon substituting in that the first Chern class $\sum t_a$ vanishes yields $$C_2(e_i) = \binom{n-2}{i-1} c_2.$$

The method of this computation extends to all symmetric polynomials: for example, it is even more straightforward to compute the Chern classes of the representation corresponding to the power sum $p_j = \sum t_a^j$ as we immediately obtain \begin{equation*} ch(p_j) = \sum \exp(t_a)^j = \sum \exp(jt_a), \end{equation*} which in turn easily yields \begin{equation} \label{power} C_i(p_j) = j^i c_i.\end{equation} In particular, for the current problem of calculating the index of $CH^2(BG)$ inside $H^4(BSL_n; \mb{Z})$, the representation corresponding to the power sum $p_{\ell}$ is certainly a representation which descends from $SL_n$ to $G$. As it satisfies $C_2(p_{\ell}) = \ell^2 c_2$, we see that the index is only one of three possibilities: $1$, $\ell$, or $\ell^2$. We claim it is in fact always $\ell$, except for the case $\ell = 2$ and $k$ odd, in which case the index is $4$. We hence should exhibit some representation $\rho$ such that $C_2(\rho)$ is not divisible by $\ell^2$, as well as estabish that for every representation, $C_2(\rho)$ is divisible by $\ell$. We first address the latter question.

Consider $C_2$ evaluated on a general basis element $\prod_{\alpha \in S} e_{\alpha}$. Indeed, it suffices to consider $C_2$ on an integral spanning set, as $C_2$ is certainly additive. Unfortunately, it is not multiplicative, but the Chern character is, so using~\ref{chern} again, we find that given two representations $\rho_1, \rho_2$, the fact that $ch(\rho_1 \rho_2) = ch(\rho_1) ch(\rho_2)$ implies that $C_2$ behaves essentially as a derivation: $$C_2(\rho_1 \rho_2) = \rank(\rho_1) C_2(\rho_2) + C_2(\rho_1) \rank(\rho_2).$$ The simplicity of this expression depended on the vanishing of $C_1$. Applying this derivation property, we arrive at the following: \begin{eqnarray*} C_2\Big(\prod_{\alpha \in S} e_{\alpha}\Big) &=& \sum_{\alpha \in S} C_2(e_{\alpha}) \prod_{\beta \in S \setminus \{\alpha\}} \rank(e_{\beta}) \\ &=& \sum_{\alpha \in S} \binom{n-2}{\alpha-1} \prod_{\beta \in S \setminus \{\alpha\}} \binom{n}{\beta}. \end{eqnarray*} We now use two elementary observations from Lucas's theorem on binomial coefficients: as $n$ is a multiple of $\ell$, if $i$ is a multiple of $\ell$, we have $$\ell \bigm| \binom{n-2}{i-1}.$$ Indeed, when expanded in base $\ell$, $n-2$ has last digit $\ell - 2$, which is smaller than the last digit $\ell - 1$ of $i-1$. On the other hand, if $i$ is not a multiple of $\ell$, then $$\ell \bigm| \binom{n}{i}$$ as now the last digits are $0$ and something nonzero, respectively. Turning to our expression for $C_2$, we see that from the $\binom{n}{\beta}$ factors, the only way to get a single term of that sum to be nonzero modulo $\ell$ is if there exists some $\alpha$ such that for all $\beta \in S \setminus \{\alpha\}$, $\beta$ is divisible by $\ell$. However, as $$\ell \bigm| \sum_{\alpha \in S} \alpha,$$ we obtain that $\alpha$ is also divisible by $\ell$, i.e. every member of $S$ is divisible by $\ell$. Hence the first factor of $\binom{n-2}{\alpha-1}$ is always divisible by $\ell$ and thus so is $C_2$.

It remains to verify that the index of $CH^2(BG)$ inside $H^4(BSL_n; \mb{Z})$ is not $\ell^2$ save for the one exception we noted above, so we simply consider a few representations, where we now consider only the case $\ell$ odd. As $e_{\ell}$ is a representation that descends to $G$, we may use $C_2(e_{\ell}) = \binom{n-2}{\ell-1}$, which by Kummer's theorem is divisible by $\ell^2$ if and only if we have at least two carries in writing $k \ell - 2$ as a sum of $\ell - 1$ and $(k-1) \ell - 1$, which happens if and only if $k \equiv 1 \pmod{\ell}.$ In that case, we instead consider $e_1e_{\ell-1}$, which has \begin{eqnarray*} C_2(e_1 e_{\ell-1}) &=& \binom{n}{1} \binom{n-2}{\ell-2} + \binom{n}{\ell-1} \\ &=& \binom{n}{\ell - 1} \frac{n \ell - \ell^2 + 2 \ell - 2}{n - 1}. \end{eqnarray*} Here the fraction is an $\ell$-adic unit, while using Kummer's theorem for $\binom{n}{\ell-1}$ reveals that we want to write, in base $\ell$, we want to write $\cdots 1 0$ as $\cdots 0 1$ plus $\ell - 1$, which clearly takes only one carry. Hence $C_2(e_1 e_{\ell-1})$ is not divisible by $\ell^2$. Finally, for the case $\ell = 2$ and $k$ odd, we could revisit the representations above and establish that $C_2$ is divisible by $4$, but we in fact already know this result from the computation of $H^4(BG, \mb{Z}) \subset H^4(BSL_n, \mb{Z})$, which in this case is indeed of index $4$. 

We have hence completely calculated the indices of $CH^2(BG) \subset H^4(BG; \mb{Z}) \subset H^4(BSL_n; \mb{Z})$, and we have a discrepancy in the indices and hence counterexamples to the integral Hodge conjecture precisely when $k$ is divisible by $\ell'$, in which case $H^4(BG; \mb{Z}) \to H^4(BSL_n; \mb{Z})$ is an isomorphism while $CH^2(BG)$ is of index $\ell$. 

\section{The Third Degree}

We turn now to a consideration of the group $Ch^3(BG)$ for $G$ still in the above generality, before focusing more carefully on the particular case of $G = SL_6 / \mu_2$. Certainly, the easy computation of~\ref{power} for $i = 3$ yields $C_3(p_{\ell}) = \ell^3$, and so the index of $\im CH^3(BG)$ in $H^6(BSL_n, \mb{Z})$ is either $1, \ell, \ell^2,$ or $\ell^3$. Let us first restrict to the case of $\ell$ odd. We first note that, because of the Riemann-Roch without denominators argument, $Ch^3 \subset CH^3$ has index at most $2$, and so here we know $\im Ch^3 = \im CH^3$ immediately. We claim that the index is always $\ell^2$, except for the exception of $\ell = 3$, $k \equiv 1 \pmod{3}$, where the index is $3^3$. Let us first show $C_3(\rho)$ is always divisible by $\ell^2$; as before, we first compute $C_3$ for the totally antisymmetric representations. We revisit the computation that led to~\ref{2}, now proceeding one degree further. We may similarly calculate, after substituting in the vanishing of the first Chern class, that $$  \frac{1}{2} C_3(e_i) = \frac{1}{2} \binom{n-1}{i-1} - \frac{3}{2} \binom{n-2}{i-2} + \binom{n-3}{i-3},$$ and so $$C_3(e_i) = \binom{n-2}{i-1} \frac{n-2i}{n-2}.$$ As before, $C_3$ also satisfies a derivation-like property, and so we have on a general element of our integral basis that \begin{equation} \label{expression} C_3 \Big( \prod_{\alpha \in S} \Big) = \sum_{\alpha \in S} \binom{n-2}{\alpha - 1} \frac{n-2\alpha}{n-2} \prod_{\beta \in S \setminus \{\alpha\}} \binom{n}{\beta}.\end{equation}

Let us first check that $C_3(e_{i\ell})$ is divisible by $\ell^2$. Note that throughout, as we assumed $\ell$ odd, we have that $n - 2$ is an $\ell$-adic unit and so does not affect $\ell$-valuation. However, as in our previous discussion of $C_2$, the binomial coefficient is certainly divisible by $\ell$, as is the numerator $n - 2i \ell$, and so we have the conclusion. It now suffices to consider multisets $S$ where all elements are prime to $\ell$; indeed, suppose we have the result for a multiset $S$. As we know the result for $S$ and $i \ell$ independently, the derivation property yields that we also know the result for the multiset given by adjoining $i \ell$ to $S$. So, suppose that $S$ has all elements coprime to $\ell$. Recall that $\binom{n}{\beta}$ is divisible by $\ell$ for $\beta$ not divisible by $\ell$, and so if we want to consider the expression~\ref{expression} modulo $\ell^2$, the only terms which can contribute are those for which at most one $\beta$ is not divisible by $\ell$. In other words, there are at most two elements not divisible by $\ell$. As the sum of all elements is divisible by $\ell$, we must have exactly two indices as opposed to merely one, so denoting these indices by $i$ and $j$, we study \begin{eqnarray*} C_3(e_i e_j) &=& \binom{n-2}{i-1} \frac{n-2i}{n-2} \binom{n}{j} + \binom{n-2}{j-1} \frac{n-2j}{n-2} \binom{n}{i} \\ &=& \frac{1}{(n-1)(n-2)} \frac{1}{n} \binom{n}{i} \binom{n}{j} \Big( i(n-i)(n-2i) + j(n-j)(n-2j) \Big) \\ &=& \frac{1}{(n-1)(n-2)} \frac{1}{n} \binom{n}{i} \binom{n}{j} (n-i-j) \Big( (i+j)n - 2(i^2 + ij + j^2) \Big). \end{eqnarray*} If $n$ is divisible by $\ell^2$, we have by Kummer's theorem that $\binom{n}{i}, \binom{n}{j}$ are already divisible by $\ell^2$, so this expression is certainly divisible by $\ell^2$; otherwise, the denominators above only detract one power of $\ell$, but $\binom{n}{i}, \binom{n}{j},$ and $n-i-j$ are all divisible by $\ell$, so the product above is divisible by $\ell^2$, as desired.

To establish that in most cases, the index is indeed $\ell^2$ rather than $\ell^3$, consider first $C_3(e_{\ell}) = \binom{n-2}{\ell - 1} \frac{n-2\ell}{n-2}$. As discussed previously with $C_2(e_{\ell})$, the binomial coefficient is not divisible by $\ell^2$ unless $k \equiv 1 \pmod{\ell}$, and similarly the fraction only introduces extra powers of $\ell$ if $k \equiv 2 \pmod{\ell}$. So we already have the result except for $k \equiv 1, 2 \pmod{\ell}$. If $k \equiv 2 \pmod{\ell}$, consider $$C_3(e_1 e_{\ell - 1}) = \frac{1}{(n-1)(n-2)} \binom{n}{\ell - 1} \Big( (k-1) \ell \Big)  \Big( \ell n  - 2 (\ell^2 - \ell + 1) \Big).$$ The last factor is coprime to $\ell$, while the penultimate factor has exactly one factor of $\ell$. By Kummer's theorem, the binomial coefficient also has only one factor of $\ell$, as there is only carry in adding $\cdots 11$ to $\ell - 1$ to obtain $\cdots 20$ in base $\ell$. Hence we have that this Chern class is exactly divisible by two powers of $\ell$. Finally, for $k \equiv 1 \pmod{\ell}$, consider \begin{eqnarray*} C_3(e_1^2 e_{\ell - 2}) &=& 2C_3(e_1) \rank(e_1) \rank(e_{\ell-2}) + C_3(e_{\ell-2}) \rank(e_1)^2 \\ &=& 2 n \binom{n}{\ell-2} + n^2 \binom{n-2}{\ell-3} \frac{n-2\ell+4}{n-2} \\ &=& \frac{n}{(n-1)(n-2)} \binom{n}{\ell - 2} \Big( \ell n^2 - 3(\ell^2 - 4\ell + 6) n + 2(\ell^3 - 6\ell^2 + 12 \ell - 6) \Big). \end{eqnarray*} Again, $n$ has exactly one factor of $\ell$, while so does the binomial coefficient by Kummer's theorem. The last factor is coprime to $\ell$ using the assumption that $\ell > 3$.

Finally, treating the case of $\ell = 3$ and $k \equiv 1 \pmod{3}$ is more involved. We leave this discussion to the appendix for the interested reader.

\section{$SL_6 / \mu_2$}

We now turn to the case $\ell = 2$ for the particular case $k = 3$ to further improve focus given the wider range of questions we consider here. We start as before by computing the index of $\im Ch^3(BG)$. In fact, we have $C_3(e_2) = \binom{4}{1} \frac{2}{4} = 2$, so the index is either $1$ or $2$. We quickly dispatch this ambiguity by explicitly calculating that for $i$ even in the relevant range, we have $C_3(e_i)$ even: indeed, $C_3(e_2) = 2, C_3(e_4) = -2$. It follows that we need only consider multisets $S$ with all indices odd, but then as the multiset must have at least two elements and $\rank(e_i) = \binom{6}{i}$ is even for $i$ odd, we see that $C_3$ is indeed always even and the index is $2$. As we perform these computations, let us define, without motivation for now, the subgroup $Ch^3_{C_2=0}(BG) \subset Ch^3(BG)$ generated by Chern classes of those virtual representations $\rho$ which satisfy $C_2(\rho) = 0$. In fact, as~\ref{power} yields that $C_2(p_2) = 4, C_3(p_2) = 8$, and we may further calculate $C_2(e_2) = 4, C_3(e_2) = 2$. As such, the virtual representation corresponding to $p_2 - e_2$ has vanishing second Chern class but $C_3(p_2 - e_2) = 6$. On the other hand, $C_2(e_4) = 4, C_3(e_4) = -2$, so $C_2(p_2 - e_4) = 0, C_3(p_2 - e_4) = 10$. Hence $\im Ch^3_{C_2=0} = \im Ch^3$ is also of index $2$ inside $H^6(BSL_6; \mb{Z})$. 

We wish now to consider not only the cycle class map $CH^3(BG) \to H^6(BG; \mb{Z})$ but also the refined cycle class map $CH^3(BG) \to MU^6(BG) \otimes_{MU_*} \mb{Z}$, ideally with a view towards Yagita's conjecture on whether this map is always an isomorphism for classifying spaces of affine algebraic groups. As usual, we only consider these groups as mapped to $H^6(BSL_n; \mb{Z})$. Moreover, as the problem is evidently $2$-local, we may further consider the interpolating cohomology theories $MU^6(BG) \to BP^6(BG) \to BP\langle 1 \rangle^6(BG) \to H^6(BG; \mb{Z})$, where the (truncated) Brown-Peterson spectra are at the prime $\ell = 2$. Here $BP$ is Brown-Peterson cohomology and $BP \langle 1 \rangle$ is one of the truncated Brown-Peterson cohomologies introduced by~\cite{Johnson-Wilson}, although the first truncation has the alternative expression $BP \langle 1 \rangle \simeq ku$, the connective cover of $K$-theory ($\ell$-locally completed). Using Wilson's main theorem in~\cite{Wilson} for $k = 6, n = 1, p = 2$, we find that $BP^6(BG) \to BP \langle 1 \rangle^6(BG)$ is surjective and hence have the same image in integral cohomology. As such, the chain of subgroups $\im Ch^3_{C_2=0}(BG) \subset \im Ch^3(BG) \subset \im CH^3(BG) \subset \im MU^6(BG) \subset \im BP^6(BG) \subset \im ku^6(BG)$ collapses somewhat, with the first inclusion and the last two inclusions all equalities.

We now turn to identifying $\im ku^6(BG)$. We use the cofibre sequence of spectra \begin{equation} \label{seq} \Sigma^2 ku \stackrel{v}{\to} ku \to H \mb{Z}\end{equation} where $v$ denotes the Bott map; we also recall that inverting the Bott element yields nonconnective K-theory, both as abstract spectra $ku[v^{-1}] \simeq K$, and more concretely when evaluating on a given CW complex $K^*(X) \simeq ku^*(X)[v^{-1}]$. Finally, we have the Atiyah-Segal completion theorem $K(BG) \simeq \Rep(G)^{\wedge}_I$, where we are completing the representation ring at the augmentation ideal; in particular, the odd $K$-theory of $BG$ vanishes. Note, for example, that these last two facts imply that any odd-degree $ku$-cohomology class is $v$-nilpotent, i.e. dies after multiplication by a sufficiently high power of the Bott class. We also need some computation of the low-degree cohomology of $BG$, which follows from the Serre spectral sequence for the fiber sequence $B\mu_2 \to BSL_6 \to BG$: $H^1(BG; \mb{Z}) = H^2(BG; \mb{Z}) = 0, H^3(BG; \mb{Z}) = \mb{Z}/2,$ and $H^4(BG; \mb{Z}) = \mb{Z}$. Although we are more interested here in even-degree $ku$-cohomology, we first note the simpler application of these facts to odd-degree $ku$-cohomology: the above computations immediately yield that the Bott maps $ku^3(BG) \stackrel{v}{\to} ku^1(BG) \stackrel{v}{\to} ku^{-1}(BG) \stackrel{v}{\to} \cdots$ are all isomorphisms, as the relevant kernels and cokernels in the long exact sequence resulting from~\ref{seq} all vanish, at least once one observes that $ku^0(BG) \to H^0(BG; \mb{Z})$ is surjective. Hence, as odd-degree $ku$-cohomology classes on $BG$ are $v$-nilpotent, we immediately have that all the above groups vanish and hence, so does $ku^5(BG)$, as it is sandwiched between two vanishing groups in the long exact sequence. Meanwhile, for the even-degree $ku$-cohomology, we claim that at least if we ignore torsion, all the Bott maps are isomorphisms, and in particular, the long exact sequence~\ref{seq} yields $$ku^6(BG) / tors \hookrightarrow ku^4(BG) \stackrel{\sim}{\to} ku^2(BG) \hookrightarrow ku^0(BG) \stackrel{\sim}{\to} ku^{-2}(BG) \stackrel{\sim}{\to} \cdots.$$ The argument proceeds by the usual observation that all possible kernels, and cokernels where relevant, vanish by the computation of integral cohomology of $BG$. Only the first injection needs further comment: the point is that any possible kernel of $ku^6(BG) \to ku^4(BG)$ is a quotient of the torsion group $H^3(BG; \mb{Z})$. Indeed, although we do not need to consider higher degrees for the present investigation, we know in general that the odd-degree integral cohomology of $BG$ is torsion as $H^*(BG; \mb{Q})$ is well-known to identify with $H^*(BT; \mb{Q})^W$ for the Weyl group action on a maximal torus, and in particular is concentrated in even degrees. Moreover, the above analysis of the Bott maps and the comparison of $ku$ to $K$ yields that, for example, $ku^4(BG)$ may be identified with a suitable subgroup of the torsion-free $\Rep(G)^{\wedge}_I$ and so $ku^6(BG)$ splits non-canonically as a further subgroup of $\Rep(G)^{\wedge}_I$ plus the torsion group $H^3(BG; \mb{Z})$. We are content to ignore this torsion subgroup as we are interested in the image of $ku^6(BG) \to H^6(BG; \mb{Z}) \to H^6(BSL_6; \mb{Z}) \simeq \mb{Z} c_3$ and a torsion summand of the source certainly admits no interesting maps to the target. Furthermore, the completion of the representation ring makes no difference to the image inside $\mb{Z} c_3$ and so we ignore it. Hence, we now wish to do three things: (i) concretely characterize the subgroup of $\Rep(G)$ which may be identified, after completion, with $ku^6(BG) / tors$ under the Bott maps, (ii) explicitly identify the map from this subgroup to $H^6(BSL_6; \mb{Z})$, and (iii) calculate the image.

We treat the first two questions by comparison of the cofibre sequence~\ref{seq} for $BG$ and $BSL_n$. Certainly, by commutativity of the squares \begin{equation*} \xymatrix{ ku^*(BG) \ar[r] \ar[d] & H^*(BG; \mb{Z}) \ar[d] \\ ku^*(BSL_n) \ar[r] & H^*(BSL_n; \mb{Z}) } \end{equation*} for $* = 4, 0$, we see that the conditions for a $G$-representation to survive to $ku^6(BG)$ are precisely the same conditions that an $SL_n$-representation to survive to $ku^6(BSL_n)$, and so we simply study the situation for $SL_n$ and only thereafter consider the subring of representations that also descend to $G$ (i.e., have even number of total boxes). However, it is easy to understand the map $ku \to H\mb{Z}$ for $BSL_n$ as the Atiyah-Hirzebruch spectral sequence computing $ku$ from $H\mb{Z}$ completely collapses due to being concentrated in even degrees. In other words, if we continue to be lazy about denoting the appropriate completions, we have $\Rep(SL_n)$ as the degree-zero part of $H^*(BSL_n; \mb{Z}) \otimes K_*$, which therefore decomposes as \begin{equation} \label{decomp} \Rep(SL_n) \simeq \mb{Z} \oplus \mb{Z} v^2 c_2 \oplus \mb{Z} v^3 c_3 \oplus \cdots.\end{equation} Moreover, the map $ku^{2i}(BSL_n) \to H^{2i}(BSL_n; \mb{Z})$ is given by projection on the $\mb{Z} v^i c_i$ factor, as we see by comparing the degenerate Atiyah-Hirzebruch spectral sequences with the cofibre sequence~\ref{seq}: when we identify $ku^{2i}(BSL_n)$ with a subgroup of $\Rep(SL_n)$ by applying the Bott map $i$ times, we precisely obtain the subgroup such that in~\ref{decomp}, the first $i$ components vanish, i.e. we obtain the subgroup of elements divisible by $v^i$, and then the map to $H\mb{Z}$ extracts the lowest component.

From this point of view, the map $\Rep(SL_n) \simeq ku^0(BSL_n) \to H^0(BSL_n; \mb{Z})$ may be identified with projection on the first factor in~\ref{decomp} (which corresponds to taking the rank of the virtual representation), so that the subgroup which survives to $ku^2$ is $\Rep(SL_n)_{\rank = 0}$. This same subgroup survives to $ku^4$, but then admits a map to $H^4(BSL_n; \mb{Z})$ given by projection onto the second factor in~\ref{decomp}. Given a virtual representation $\rho$ in the subgroup $\Rep(SL_n)_{\rank = 0}$, we wish to understand what this map $c_2$ really means. We do so by thinking more carefully about how the Atiyah-Hirzebruch spectral sequence is constructed: namely, by taking a Postnikov decomposition of the target. From this point of view, we start with our virtual representation $\rho$, which we interpret as a map $BSL_n \to \mb{Z} \times BU \simeq ku$, where how you write the target differs depending on whether you prefer to equivalently think of the infinite loop space or the connective spectrum. We will use the dual Whitehead tower, so consider the following diagram: \begin{equation*} \xymatrix{ & \vdots \ar[d] & \\ & BU^{(6)} \ar[d] \ar[r] & K(\mb{Z}, 6) \\ & BU^{(4)} \ar[d] \ar[r] & K(\mb{Z}, 4) \\ & BU^{(2)} \ar[d] \ar[r] & K(\mb{Z}, 2) \\ BSL_n \ar[r] \ar@{-->}[ur] \ar@{-->}[uur] \ar@{-->}[uuur] & \mb{Z} \times BU \ar[r] & K(\mb{Z}, 0). } \end{equation*} Here $BU^{(n)}$ represents the $n$-connective ($(n-1)$-connected) cover of $BU$, so for example $BU^{(2)}$ is simply $BU$ itself; the arrows to the Eilenberg-MacLane spaces represent their lowest degree cohomology, or equivalently their lowest degree homotopy by Hurewicz, which is why we see the homotopy groups of the original target splayed out on the right-hand column, as we expected from the Atiyah-Hirzebruch $E_2$-page. Now, the way the Atiyah-Hirzebruch spectral sequence works in the simplest case when it completely degenerates is that we first have a contribution from maps to the bottom Eilenberg-MacLane space which represents projection to the first factor in a decomposition such as~\ref{decomp}. If this projection vanishes, the map lifts to the next, more highly-connected cover in the Whitehead tower and we pick up a contribution from maps to the next Eilenberg-MacLane space, and so on. In other words, given that $\rank \rho = 0$ and so our map lifts to the zero component of $\mb{Z} \times BU$, the projection on the next factor is precisely given by composing with $c_2: BU \to K(\mb{Z}, 2)$. This composition is precisely what we referred to as $C_2(\rho)$, and so we also have our answer that a representation lifts to $ku^6(BSL_n)$ is this projection vanishes, i.e. if both $\rank(\rho)$ and $C_2(\rho) = 0$. These conditions characterize the subgroup of $\Rep(SL_n)$ that we care about, and so we have our answer to the first of our three missions.

We are in prime condition to embark on the next desideratum: identifying the map to $H^6(BSL_n; \mb{Z})$ on this subgroup is exactly the problem of understanding the next projection in~\ref{decomp} to $\mb{Z} v^3 c_3$, so we wish to understand the map $f: BU^{(6)} \to K(\mb{Z}, 6)$. In fact, we already have such a map, given by the composition $g: BU^{(6)} \to BU \stackrel{c_3}{\to} K(\mb{Z}, 6)$, and we already understand this latter map applied to a representation, as it is precisely what we denoted $C_3(\rho)$. To understand how the two maps compare, we know that $g$ must be some integral multiple of $f$, as $f$ classifies a generator of $H^6(BU^{(6)}; \mb{Z})$, which by Hurewicz we know is dual to $\pi_6(BU^{(6)} \simeq \pi_6 BU \simeq \mb{Z}$. In fact, as the latter map classifies a generator for $H^6(BU; \mb{Z}) \simeq \mb{Z}$, this argument informs us that the multiple we seek is the same as the index of the Hurewicz map $\pi_6(BU) \to H_6(BU; \mb{Z})$. We can calculate this index by evaluating the third Chern class on the cube of the Bott element, but as the Chern character evaluated on a power of the Bott element is simply $Ch(v^n) = Ch(v)^n = (e^t - 1)^n = t^n + O(t^{n+1})$, where $t$ is the first Chern class of the tautological line bundle on $S^2$, we see that the index we seek is precisely the denominator of the coefficient of $c_3$ when expanding the Chern character in terms of Chern classes. By the same numerical manipulations as goes into Riemann-Roch without denominators, we hence see that in general, the composite $BU^{(2i)} \to BU \stackrel{c_i} \to K(\mb{Z}, 2i)$ is $(-1)^{i-1} (i-1)!$ times the generator $BU^{(2i)} \to K(\mb{Z}, 2i)$, and so we have now concluded with the second task: the projection to $\mb{Z} v^3 c_3$ on a representation $\rho$ satisfying the constraints $\rank \rho = 0, C_2(\rho) = 0$ is given by $\frac{1}{2} C_3(\rho)$. In particular, representations satisfying the given constraints always have even $C_3$. Indeed, this careful consideration of the Postnikov or Whitehead tower of $BU$ is an elementary version of, for example, Adams's approach to finding a much larger class of such divisibility relations in~\ref{Adams}.

We now turn to the third task of actually identifying the image of this map. Ideally our earlie definition of $Ch^3_{C_2 = 0}(BG)$ is now motivated: we are precisely interested in the representation with vanishing second Chern class in order to evaluate their third Chern class and divide it by two. Note that the condition on the rank vanishing is fairly superfluous, as we can always correct by subtracting a suitable number of trivial representations. Hence we now have the equality of subgroups $\im ku^6(BG) = \frac{1}{2} \im Ch^3_{C_2 = 0}(BG)$, so in the case at hand, the index of this subgroup is $1$ and $ku^6(BG)$ surjects onto $H^6(BSL_6; \mb{Z})$. We hence have the chain of subgroups $2 \mb{Z} = \im Ch^3(BG) \subset \im CH^3(BG) \subset \im MU^6(BG) = \mb{Z}$. It remains to address the issue of whether $\im CH^3(BG)$ is indeed $\mb{Z}$ or $2 \mb{Z}$. In the former case, we have another nice example of a Chow group of a classifying space not generated by Chern classes of polynomial representations. In the latter case, we have a counterexample to Yagita's conjecture. Hence, the question is certainly of pressing interest; however, we cannot treat it with the tools of representation theory and algebraic topology we brought to bear in this paper, as we need some further legitimate geometry. In particular, the problem looks amenable to the stratification method used to great effect by Vistoli and Vistoli-Molina Rojas in~\cite{Vistoli} and~\cite{Vistoli-MR}; indeed, we expect in future work to analyze the orbit stratification of the second symmetric power of the standard representation.

\section{Appendix: $CH^3(BSL_{9j+3}/\mu_3)$}

We continue the discussion of section $4$ for the case $\ell = 3$, $k \equiv 1 \pmod{3}$. We wish to show $C_3(\rho)$ is always divisible by $\ell^3$. We start by straightforwardly establishing the claim for $\rho = e_{i \ell}$: in the expression $C_3(e_{i \ell}) = \binom{n-2}{i \ell - 1} \frac{n - 2 i \ell}{n-2}$, the fraction has a power of $\ell$ while the binomial coefficient has at least two powers of $\ell$ by Kummer's theorem as discussed in section $3$. As before, it now suffices to establish our claim for multisets $S$ all of whose elements are prime to $\ell$. Also as before, as $\ell \bigm| \binom{n}{\beta}$ for $\beta$ coprime to $\ell$, we can have at most three indices prime to $\ell$ before the binomial coefficients make the expression automatically divisible by $\ell^3$, and so subject to our assumption on $S$, it has either two or three indices. If we have only two indices in $S$, we copy the result of our previous computation: $$C_3(e_i e_j) = \frac{1}{(n-1)(n-2)} \frac{1}{n} \binom{n}{i} \binom{n}{j} (n-i-j) \Big((i + j)n - 2(i^2 + ij + j^2)\Big).$$ As $\val_{\ell} n = 1$, it suffices to establish $\ell^4 \bigm| \binom{n}{i} \binom{n}{j} (n-i-j)$. We certainly have that all of these factors are divisible by $\ell$, but we claim that in fact one of them must be divisible by $\ell^2$. Indeed, suppose otherwise. If we analyze the binomial coefficient $\binom{n}{i}$ by Kummer's theorem, note that $n$ expanded in base $\ell$ ends with $\cdots 1 0$. The only way to write this as a sum of two numbers prime to $\ell$ and have precisely one carry is if the penultimate digit of $i$ is $0$, or in other words $i$ reduced modulo $\ell^2$ is smaller than $\ell$. As $j$ satisfies the same condition, we see that $i + j \equiv \ell \pmod{\ell^2}$ and hence $\ell^2 \bigm| n-i-j$, establishing the claim.

On the other hand, if $S$ has three indices $i, j,$ and $m$, then a similar computation yields $$C_3(e_ie_je_m) = \frac{1}{(n-1)(n-2)} \frac{1}{n} \binom{n}{i} \binom{n}{j} \binom{n}{m} \Big( i(n-i)(n-2i) + j(n-j)(n-2j) + m(n-m)(n-2m) \Big).$$ The fraction removes one power of $\ell$ while the binomial coefficients are responsible for at least three. We use that $\ell = 3$ to claim that the last factor is also divisible by $\ell$: working modulo $\ell$, we have \begin{eqnarray*} i(n-i)(n-2i) + j(n-j)(n-2j) &\equiv& 2(i^3 + j^3 + m^3) \pmod{\ell} \\ &\equiv& 2(i^3 + j^3 - (i+j)^3) \pmod{\ell} \\ &\equiv& -6ijm \pmod{\ell},\end{eqnarray*} whence the result. Based on these computations, we suspect in general that $CH^p(BSL_{kp} / \mu_p)$ for $k \equiv 1 \pmod{p}$ is a particularly rich place to seek counterexamples to the integral Hodge conjecture; we hope to later investigate this suspicion and the properties of the putative cycles in question as well. 

\bibliographystyle{alpha}
\bibliography{ref}

\begin{thebibliography}{MRV06}

\bibitem[AH62]{Atiyah-Hirzebruch}
M.~F. Atiyah and F.~Hirzebruch.
\newblock Analytic cycles on complex manifolds.
\newblock {\em Topology}, 1:25--45, 1962.

\bibitem[Ant15]{Antieau}
Ben Antieau.
\newblock On the integral tate conjecture for finite fields and representation
  theory.
\newblock \url{http://arxiv.org/abs/1504.04879/}, 2015.

\bibitem[AW14]{Antieau-Williams}
Benjamin Antieau and Ben Williams.
\newblock The topological period-index problem over 6-complexes.
\newblock {\em J. Topol.}, 7(3):617--640, 2014.

\bibitem[JW73]{Johnson-Wilson}
David~Copeland Johnson and W.~Stephen Wilson.
\newblock Projective dimension and {B}rown-{P}eterson homology.
\newblock {\em Topology}, 12:327--353, 1973.

\bibitem[Kam15]{Kameko}
Masaki Kameko.
\newblock On the integral {T}ate conjecture over finite fields.
\newblock {\em Math. Proc. Cambridge Philos. Soc.}, 158(3):531--546, 2015.

\bibitem[Kol92]{KollarHodge}
J{\'a}nos Koll{\'a}r.
\newblock Trento examples.
\newblock {\em Classification of irregular varieties (Trento, 1990)},
  1515:136--139, 1992.

\bibitem[MRV06]{Vistoli-MR}
Luis~Alberto Molina~Rojas and Angelo Vistoli.
\newblock On the {C}how rings of classifying spaces for classical groups.
\newblock {\em Rend. Sem. Mat. Univ. Padova}, 116:271--298, 2006.

\bibitem[PY14]{Pirutka-Yagita}
Alena Pirutka and Nobuaki Yagita.
\newblock Note on the counterexamples for the integral tate conjecture over
  finite fields.
\newblock \url{http://arxiv.org/abs/1401.1620/}, 2014.

\bibitem[Tot97]{TotaroCycleClass}
Burt Totaro.
\newblock Torsion algebraic cycles and complex cobordism.
\newblock {\em J. Amer. Math. Soc.}, 10(2):467--493, 1997.

\bibitem[Tot99]{TotaroClassifying}
Burt Totaro.
\newblock The {C}how ring of a classifying space.
\newblock In {\em Algebraic {$K$}-theory ({S}eattle, {WA}, 1997)}, volume~67 of
  {\em Proc. Sympos. Pure Math.}, pages 249--281. Amer. Math. Soc., Providence,
  RI, 1999.

\bibitem[Tot14]{TotaroBook}
Burt Totaro.
\newblock {\em Group cohomology and algebraic cycles}, volume 204 of {\em
  Cambridge Tracts in Mathematics}.
\newblock Cambridge University Press, Cambridge, 2014.

\bibitem[Vis07]{Vistoli}
Angelo Vistoli.
\newblock On the cohomology and the {C}how ring of the classifying space of
  {${\rm PGL}_p$}.
\newblock {\em J. Reine Angew. Math.}, 610:181--227, 2007.

\bibitem[Wil75]{Wilson}
W.~Stephen Wilson.
\newblock The {$\Omega $}-spectrum for {B}rown-{P}eterson cohomology. {II}.
\newblock {\em Amer. J. Math.}, 97:101--123, 1975.

\bibitem[Yag10]{Yagita}
Nobuaki Yagita.
\newblock Coniveau filtration of cohomology of groups.
\newblock {\em Proc. Lond. Math. Soc. (3)}, 101(1):179--206, 2010.

\end{thebibliography}
\end{document}